\documentclass{amsart}
\usepackage[T1]{fontenc}
\usepackage{amssymb}
\usepackage[hyperref]{degt}

\def\2{\color{red}}

\def\minitab#1{\vcenter\bgroup\rm
 \def\minitabskip{#1}%
 \def\-##1{\setbox0\hbox{$00$}\hbox to\wd0{\hss$##1$\hss}}%
 \def\.{\cdot}
 \let\\\cr
 \halign\bgroup\strut\hss$##$&&#1\hss$##$\hss\cr}
\def\endminitab{\crcr\egroup\egroup}

\def\trule{\noalign{\vspace{1pt}\hrule\vspace{2pt}}}
\def\HHH{\setbox0\hbox{$(0)$}\hbox to\wd0{\hss$\HH$\hss}}
\def\counts#1{\multispan2\strut\hss#1:\ignorespaces}
\def\configurations{\counts{configurations}}
\def\clist#1{\omit\minitabskip$#1$\hss}
\def\tabconfig#1{\hbox to12pt{\hss$#1$\!\!\hss}}

\def\intr(#1x#2){(#1\times#2)}
\def\GRAPH{%
 \def\null{}%
 \def\sep##1{\ssep\let\ssep=##1 }
 \def\AA[##1]{\sep\ \tA_{##1}}%
 \def\DD[##1]{\sep\ \tD_{##1}}%
 \def\A[##1]{\sep\ \bA_{##1}}%
 \def\(##1x##2){\sep\null\intr(##1x##2)}%
 \let\ssep\null
}

\def\rminus{\mathbin{{\sminus}\!\!\!{\sminus}}}
\def\HH{\Sigma}
\def\girth{\operatorname{girth}}

\let\PL\relax
\def\PIC{%
 \setlength\unitlength{1.2pt}%
 \def\b{\circle*4}%
 \def\w{\circle4}%
 \def\c{\circle6}%
 \def\(##1,##2){\put(##1,##2)}%
 \def\PLL{\smash{\llap{\PL\ }}}
 \def\nono##1{\hbox to0pt{\hss$\scriptstyle##1$\hss}}
 \def\no(##1,##2)##3{\put(##1,##2){\nono{##3}}}
 \vcenter\bgroup\hbox\bgroup\begin{picture}
}
\def\endPIC{\end{picture}\egroup\egroup}

\def\Fn{\operatorname{Fn}}
\def\sat{\operatorname{sat}}
\def\hp{\operatorname{hp}}
\def\Fano{\Cal{F}}
\def\bL{\bold{L}}
\def\CK{\Cal K}
\let\att\varsigma
\def\Hnum{\#_\HH}
\def\mult{\operatorname{mult}}
\def\val{\operatorname{val}}

\def\fS{\frak{S}}

\def\subsec{\par\medskip\noindent\refstepcounter{subsubsection}%
 \underline{\thesubsubsection.}\enspace\ignorespaces}

\title[Split hyperplane sections]
{Split hyperplane sections\\ on smooth polarized $K3$-surfaces}

\author{Alex Degtyarev}

\address{%
Department of Mathematics\\
Bilkent University\\
06800 Ankara, TURKEY}

\email{
degt@fen.bilkent.edu.tr}

\thanks{%
The author was partially supported by the T\"{U}B\DOTaccent{I}TAK grant 123F111%
}



\begin{document}

\begin{abstract}
We study hyperplane sections of smooth polarized $K3$-surfaces that split
into unions of lines. We describe the dual adjacency graphs of such sections
and find
sharp upper bounds on their number. In most cases (starting from degree $6$),
we also compute all
configurations of such sections.
\end{abstract}

\maketitle

\section{Introduction}\label{S.intro}

Virtually every known proof (\cf. \cite{Segre,rams.schuett,DIS})
of Segre's celebrated theorem~\cite{Segre} on at
most $64$ lines on a smooth quartic surface $X\subset\Cp3$
contains, as an essential part,
the following argument: consider a plane section of~$X$ split into four lines
$\ell_1,\ldots,\ell_4$; each of $\ell_i$ intersects at most $18$ lines and
every line not in the plane intersects exactly one of~$\ell_i$;
hence, the number of lines is
at most $4\cdot(18-3)+4=64$. In other words, proofs make an essential use of
the plane sections of~$X$ that split into lines; in the dual adjacency graph
of lines on~$X$, the so-called \emph{Fano graph} $\Fn X$, they are
fragments isomorphic to the complete graph $K(4)$. (The proofs differ in their
approach to the configurations of lines that do \emph{not} contain a
$K(4)$-fragment and in their treatment of the very special lines that
intersects more that~$18$ others: the so-called \emph{lines of the
second kind} may have valency up to~$20$, \cf. \cite{rams.schuett:second.kind}.)

Likewise, many papers dealing with lines on sextic $K3$-surfaces
$X\subset\Cp4$ make an extensive use of the split hyperplane sections of~$X$.
Recall that $X=Q_2\cap Q_3$ is a regular complete intersection of a quadric
and a cubic. A general hyperplane~$H$ intersects~$Q_2$ in a quadric
$\Cp1\times\Cp1\subset\Cp3$, on which $Q_3$ cuts a curve~$C$ of bi-degree
$(3,3)$. If $C$ splits into lines, it must consist of six generatrices, three
from each of the two rulings, resulting in the complete
bipartite graph $K(3,3)\subset\Fn X$.

These observations lead us to the following natural question.

\problem[I. Dolgachev~\cite{DDK}]\label{problem}
What is the maximal number of
\roster
\item\label{i.4}
$K(4)$-fragments in the Fano graph of a smooth quartic $X\subset\Cp3$?
\item\label{i.6}
$K(3,3)$-fragments in the Fano graph of a smooth $K3$-sextic $X\subset\Cp4$?
\endroster
Certainly, similar questions can be asked for $2d$-polarized $K3$-surfaces
$X\subset\Cp{d+1}$.
\endproblem

To begin with, observe that a hyperplane section of a surface of degree~$2d$
is a curve of degree~$2d$ and, should it split into lines, there must be
$2d$ of them. On the other hand, due to~\cite{degt:lines}, the
maximal number of
lines on a $2d$-polarized $K3$-surface tends to decrease, oscillating between
$24$ and~$21$ when $2d\gg0$. We show that split hyperplane sections
(for short, \emph{$h$-fragments})
do not exist
if $2d>28$.

Next, before addressing \autoref{problem} in its general form, one needs to
understand the combinatorial structure of $h$-fragments as abstract graphs.
The straightforward geometric description as in degrees~$4$ and~$6$ seems no
longer available and, to our surprise, we discover that,
more often than not, this graph is not unique. Even in degree $2d=6$ the
quadric surface
$H\cap Q_2$ may split into a pair of planes and, in addition to
$K(3,3)$, we obtain the $3$-prism graph (see \autoref{fig.h=6} below).
Even though this intersection looks degenerate,
the two families have the same codimension~$5$ in the space of all
smooth sextic
$K3$-surfaces.

The following theorem is the principal result of this paper.

\theorem[see \autoref{s.idea}]\label{th.main}
The maximal number {\rm\texttt{max\,\#}} of split hyperplane sections
on a smooth $2d$-polarized $K3$-surface is given by the following table
\[*
\minitab\quad
h^2=2d          & 2& 4& 6  & 8 &10 &12&14&16&18&20&22&24&28\\\trule
\text{graphs}   & 1& 1& 2  & 3 & 6 & 9& 8& 8& 5& 3& 1& 1& 1\\
\max\#          &72&72&76  &80 &16 &90&12&24& 3& 4& 1& 1& 1\\
\text{$h$-configs}  & ?& ?&9235&860&171&44&21&12& 6& 3& 1& 1& 1\\
\endminitab
\]
We list also the number {\rm\texttt{graphs}} of combinatorial types of
$h$-fragments and, whenever known, the number {\rm\texttt{$h$-configs}} of
$h$-configurations
\rom(see \autoref{conv.config} below\rom).
\endtheorem

\autoref{th.main} is mostly proved by brute force, listing all configurations
spanned by $h$-fragments. (As a result, this gives one a way to answer a
number of further
questions about the geometry of such configurations.) We outline the principal
ideas in \autoref{S.split}, upon which, in \autoref{S.4}--\autoref{S.higher}, on
a case-by-case basis, we give the necessary details and provide some extra
information, such as the structure of the $h$-fragments, bounds itemized by the
types of the fragments, extremal configurations, \etc.
Finally, in \autoref{S.hyperelliptic} we discuss briefly the toy example of
hyperelliptic $K3$-surfaces.

\subsection{Acknowledgements}
I am grateful to Igor Dolgachev (who also suggested the original question) and
Shigeyuki Kond\=o for a number of fruitful discussions of this project at its
early stages.

\section{Split hyperplane sections}\label{S.split}

In this section, we recall the arithmetical restatement of the problem
(see \autoref{s.prelim}), classify the $h$-fragments
on birational models (see \autoref{s.fragments}),
discuss pairwise intersections of $h$-fragments
(see \autoref{s.intersections}), and outline
the idea of the computation leading to the proof of \autoref{th.main}
(see \autoref{s.idea}).

\subsection{Preliminaries \noaux{(see \cite{DIS,degt:lines,degt.Rams:octics})}}\label{s.prelim}
Given a graph~$\Gamma$ and an even integer $2d>0$, we denote
\[*
\Fano_{2d}(\Gamma):=(\Z\Gamma+\Z h)/\!\ker,\quad h^2=2d,
\]
where, for vertices $u,v\in\Gamma$, we let $u^2=-2$,
$u\cdot h=1$, and $u\cdot v=1$
(resp.\ $0$) if $u$ and~$v$ are (resp.\ are not) connected by an edge. More
generally, if $\Gamma$ is a multigraph, $u\cdot v$ is the multiplicity of the
edge $[u,v]$. We consider only those graph~$\Gamma$ for which
$\Fano_{2d}(\Gamma)$ is hyperbolic, $\Gs_+\Fano_{2d}(\Gamma)=1$.
This implies that $\Gamma$ is either
\roster*
\item
\emph{elliptic}, $\Gs_+=\Gs_0=0$ (a union of elliptic Dynkin diagrams), or
\item
\emph{parabolic}, $\Gs_+=0$, $\Gs_0>0$ (a union of Dynkin diagrams, elliptic
or affine, of which at least one is affine), or
\item
\emph{hyperbolic}, $\Gs_+=1$ (we make no assumption about $\Gs_0$).
\endroster
In the latter case, the
property that $\Gs_+\Fano_{2d}(\Gamma)=1$
depends on the degree~$2d$, \cf. the \emph{intrinsic polarization}
in~\cite{degt:lines}.

Conversely, given a $2d$-polarized hyperbolic lattice $N\ni h$, its
\emph{Fano graph} is
\[*
\Fn(N,h):=\bigl\{u\in N\bigm|u^2=-2,\ u\cdot h=1\bigr\},
\]
where two vertices $u,v\in\Fn(N,h)$ are connected by an edge of multiplicity
$u\cdot v$. If $N=\NS(X)$ for a smooth $2d$-polarized $K3$-surface
$X\subset\Cp{d+1}$ and $h\in N$ is the class of hyperplane section,
then, by the Riemann--Roch theorem, $\Fn(N,h)=\Fn X$ is
indeed the graph of lines on~$X$.

A $2d$-polarized hyperbolic lattice $N\ni h$ is called \emph{$m$-admissible},
$m=1,2,3$, if there is no vector $e\in N$ such that either
\roster*
\item
$e^2=-2$ and $e\cdot h=0$ (\emph{exceptional divisor}) or
\item
$e^2=0$ and $e\cdot h=r$ (\emph{$r$-isotropic vector}), $\ls|r|\le m$.
\endroster
An $m$-admissible lattice is called $m$-geometric if it admits a primitive
isometry to the \emph{$K3$-lattice} $\bL:=2\bE_8\oplus3\bU$. It is called
\emph{$m$-subgeometric} if it has an $m$-geometric (in particular,
$m$-admissible) finite index extension.
Recall, see \cite{Nikulin:forms}, that finite index extensions
$\tilde{N}\supset N$ are in a bijection with isotropic subgroups
$\CK:=\tilde{N}/N\subset\discr N$; those corresponding to
$m$-geometric extensions are called \emph{$m$-geometric kernels} of~$N$.

\remark\label{rem.admissible}
The $1$-admissibility of $N\ni h$ suffices to assert that $u\cdot v\ge0$ for
any pair of vertices $u,v\in\Fn(N,h)$, see~\cite{degt.Rams:octics}.
If $2d\ge4$ and $N\ni h$ is
$2$-admissible, one also has $u\cdot v\le1$ for any pair of vertices,
so that $\Fn(N,h)$ is a simple graph.
\endremark

With the degree $2d$ fixed, a graph $\Gamma$ is called
$m$-admissible/$m$-subgeometric if so is the lattice
$\Fano_{2d}(\Gamma)\ni h$. Both properties are hereditary. If
$\Gamma$ is $m$-admissible, we can speak about the
\emph{saturations}
\[*
\sat\Gamma:=\Fn\Fano_{2d}(\Gamma),\qquad
\sat(\Gamma,\CK):=\Fn\Fano_{2d}(\Gamma,\CK),
\]
where $\CK$ is an $m$-geometric kernel of $\Fano_{2d}(\Gamma)$ and
$\Fano_{2d}(\Gamma,\CK)$ is the corresponding finite index extension.
An $m$-subgeometric graph~$\Gamma$ is $m$-geometric if
$\Gamma=\sat(\Gamma,\CK)$ for an
appropriate $m$-geometric kernel~$\CK$.
Our proof of \autoref{th.main} is based on the following statement, which is
an immediate consequence of the general theory of $K3$-surfaces (the global
Torelli theorem~\cite{Pjatecki-Shapiro.Shafarevich}, surjectivity of the
period map~\cite{Kulikov:periods}, description of projective
models~\cite{Saint-Donat}, fine moduli space~\cite{Beauville:moduli}, \etc.)

\theorem[\cf. \cite{DIS,degt:lines,degt.Rams:octics}]\label{th.graphs}
A graph~$\Gamma$ is isomorphic to the
Fano graph of a
smooth $2d$-polarized $K3$-surface $X\subset\Cp{d+1}$, $2d\ge4$, if and only
if $\Gamma$ is $2$-geometric.
\done
\endtheorem

In the rest of the paper, the parameter $m$ is usually understood (typically
$m=2$ as in \autoref{th.graphs}) or not important (provided that it is
consistent throughout a statement) and we often omit it from the terminology.

\remark
In this paper, graphs represent lines on surfaces and, at the same time,
their vertices span N\'eron--Severi lattices. For this reason, we freely mix
graph theoretic,
geometric, and arithmetic terminology. Thus, for two vertices $u,v\in\Gamma$,
the terms
``are adjacent'', ``intersect'', and ``$u\cdot v=1$'' are synonyms, and so
are the terms ``are disjoint/skew'' and ``$u\cdot v=0$.''
Vertices are often referred to as lines.
\endremark

\subsection{Classification of $h$-fragments}\label{s.fragments}
Let $X\subset\Cp{d+1}$ be a smooth $2d$-polarized $K3$-surface. By
definition, an induced subgraph $\HH\subset\Fn X$
represents a split hyperplane
section (from now on, an \emph{$h$-fragment}) if $h=\sum_{v\in\HH}v$ in
$\NS(X)$.
Using \autoref{th.graphs}, we conclude that an abstract
graph~$\HH$ can serve as an $h$-fragment in degree~$2d$ if and only if
\[
\ls|\HH|=2d,\qquad
\text{$\HH$ is $3$-regular},\qquad
\text{$\HH$ is $2$-subgeometric}.
\label{eq.h-fragment}
\]
Indeed, the last condition means that $\HH$ does appear in the Fano graph
of a $K3$-surface, whereas the former two imply that the vector
$u:=h-\sum_{v\in\HH}v$ lies in the radical of the hyperbolic lattice
$\Z\HH+\Z h$; it follows that $u$ must vanish in $N/\!\ker N$ any
hyperbolic overlattice  $N\supset(\Z\HH+\Z h)$.

A \emph{bouquet} of $h$-fragments (at a line $\ell\in\Fn X$) is a union of
$h$-fragments containing~$\ell$. The \emph{multiplicity} $\mult\ell$ is the
number of $h$-fragments in the maximal bouquet at~$\ell$, \ie, the total
number of $h$-fragments $\HH\subset\Fn X$ containing~$\ell$.
The number of $h$-fragments in $\Fn X$ is denoted by $\Hnum(X)$.

\remark\label{rem.universal}
The argument above shows also that the property of~$\HH$ to be an
$h$-fragment does not depend on its embedding to $\Fn X$: whenever $\Fn X$
contains an induced subgraph isomorphic to~$\HH$, it represents a split
hyperplane section.

Furthermore, whenever $\Fn X$ contains $\HH\sminus v$, $v\in\HH$, it
also contains the missing vertex~$v$, constituting, together with
$\HH\sminus v$, an $h$-fragment.
\endremark

A $3$-regular graph cannot be elliptic or parabolic, and we employ the
taxonomy of hyperbolic graphs developed
in~\cite{degt:lines,degt.Rams:octics}. Order affine Dynkin diagrams according
to their Milnor number, followed by $\bA_\mu<\bD_\mu<\bE_\mu$. A hyperbolic
graph~$\Gamma$ contains an affine Dynkin diagram as an induce subgraph; we
fix a minimal one $\Phi\subset\Gamma$ and call~$\Gamma$ a $\Phi$-graph. This
choice splits $\Gamma$ into two subgraphs
\[*
\Pi:=\Phi\cup\bigl\{v\in\Gamma\bigm|v\cdot\Phi=0\bigr\},\qquad
\sec\Phi:=\Gamma\sminus\Pi.
\]
The connected components of the parabolic subgraph~$\Pi$ are called
\emph{fibers}, whereas the vertices of $\sec\Phi$ are called \emph{sections}.
(This terminology is motivated by elliptic pencils on $K3$-surfaces,
see~\cite{degt:lines,degt.Rams:octics} for details.)
The admissibility assumption imposes certain restrictions, \eg,
\[*
s\cdot\kappa_\Phi=s\cdot\kappa_{\Phi'}\ge s\cdot\Phi''
\]
for any section~$s$, parabolic fiber~$\Phi'$, and elliptic fiber~$\Phi''$.
(Here, $\kappa_\Phi$ stands for the fundamental cycle of the affine Dynkin
diagram~$\Phi$;
we have $s\cdot\kappa_\Phi=1$ if $s$ is a simple section.)
The assumption that $\Phi$ is minimal imposes further restrictions on the
number of sections and their pairwise intersections,
\cf.~\cite[Fig.\ 1]{degt:lines}. Finally, the $3$-regularity of~$\Gamma$
leaves us with but the following types of fibers
\roster*
\item
parabolic~$\tA_n$ and elliptic $\bA_\mu$, $\mu<n$, if $\Phi\cong\tA_n$,
$n=2,3,4,5$, or
\item
parabolic~$\tD_5$, $\tA_7$ and elliptic $\bD_4$, $\bA_\mu$, $\mu<7$,
if $\Phi\cong\tD_5$,
\endroster
asserts that there are at most three parabolic fibers,
and dictates a precise set $\sec\Phi$ of sections, \viz. the one that makes all
vertices of~$\Phi$ trivalent. (For the moment, we assume  all sections
simple.) More precisely, if $\Phi=\{a_1,\ldots,a_{n+1}\}\cong\tA_n$, there is a
unique section~$s_i$ intersecting~$a_i$ and disjoint from~$a_j$, $j\ne i$; it
is this numbering of the sections that is used in \autoref{conv.h-fragments}
below. If $\Phi=\{a_1,\ldots,a_4,b_1,b_2\}\cong\bD_5$, there are two
sections $s_{2i-1},s_{2i}$ intersecting each of the monovalent vertices~$a_i$
and disjoint from all other vertices.

It remains to run the algorithm of~\cite{degt:lines,degt.Rams:octics} and
list all graphs. Note that, since we consider smooth surfaces,
\cite[Warning~3.6]{degt.Rams:octics} does not apply: we can freely rule out all
intermediate graphs containing an affine Dynkin diagram $\Phi'<\Phi$.
We use \GAP~\cite{GAP4.13} and, most notably, its \texttt{digraph} package.
The
results of this computation are listed in the subsequent sections, separately
for each degree.

\convention\label{conv.h-fragments}
Most graphs found do not have ``standard names'', and we use an encoding
arising from our classification. Namely, we represent a graph as
\[*
\Phi_1(\sec_1)\,\Phi_2(\sec_2)\,\ldots\,\<\text{intersections}\>.
\]
Here, $\Phi_1:=\Phi,\Phi_2,\ldots$ are the fibers; their monovalent vertices
are numbered in a certain standard way. Sections are also numbered, and
$\sec_i=\sec_{i1};\sec_{i2};\ldots$, with $\sec_{ij}$ listing the sections
intersecting vertex $a_{ij}\in\Phi_i$. The optional part
\[*
\<\text{intersections}\>=\intr(n_1xn_2)\intr(n_3xn_4)\ldots
\]
lists all pairwise intersections of the sections: \eg, $s_1\cdot s_2=1$
if and only if $\intr(1x2)$ is
present on this list.

Within each section, we number the $h$-fragments $\HH_n$ according to the
list found at the beginning of this section.
In the tables where the graphs are listed, the first column ``$\HH$'' refers
to the list and the second one is ``$(r,g,s)$'', where
\roster*
\item
$r:=\rank\Z\HH=\rank\Fano_{2d}(\HH)$
or, equivalently, the codimension of the respective
stratum in the space of all smooth $2d$-polarized $K3$-surfaces,
\item
the girth $g:=\girth\HH$, and
\item
the size $s:=\ls|\Aut\HH|$ of the automorphism group (computed by
\texttt{digraph}).
\endroster
The other columns are degree specific and explained in the respective
sections.
\endconvention

It remains to consider graphs with multisections, \ie, sections~$s$ with
$s\cdot\kappa_\Phi>1$. This manual argument is also a good example of the
algorithm used above.  By~\cite{degt:lines},
$\Phi$ must be $\tA_2$
or~$\tA_3$ (or $\tD_4$, but the latter
has more than trivalent vertices).

If $\Phi=\{a_1,a_2,a_3\}\cong\tA_2$, a multisection may (and does,
see \autoref{rem.universal}) exist only if
$2d=4$, see \cite[Lemma~8.3]{degt:lines}. This section is necessarily triple,
adjacent to all three vertices of~$\Phi$,
and we obtain the $h$-fragment $K(4)$.

If $\Phi=\{a_1,a_2,a_3,a_4\}\cong\tA_3$, a bisection~$s_{13}$ may exist only if
$2d=4$, $6$, or~$8$, see \cite[Lemma~7.4]{degt:lines}. If $2d=4$, the graph
would have too many (at least~$5$) vertices. If $2d=6$, there also is another
bisection~$s_{24}$ with $s_{13}\cdot s_{24}=1$, see
\cite[Lemma~7.6]{degt:lines} or \autoref{rem.universal},
and we obtain the $h$-fragment $K(3,3)$.

Finally, if $2d=8$, then $s_{13}$ is the only bisection, which is disjoint from
all other sections, \cf.
\cite[\S\,5.3]{degt:lines}. By the $3$-regularity, there should be
two other sections $s_2,s_4$ intersecting $a_2,a_4$, respectively.
Since we need $\ls|\Gamma|=8$, this leaves but a single extra fiber
$\{a_1'\}\cong\bA_1$, and we must have
$s_{13}\cdot a_1'=s_2\cdot a_1'=s_4\cdot a_1'=1$.
By the $3$-regularity
again, we conclude that $s_2\cdot s_4=1$,
so that $\{a_1',s_2,s_4\}\cong\tA_2$ and $\Gamma$ is an $\tA_2$-graph
(without multisections): it is graph~\iref{h=8,g=1} in \autoref{S.8}
and \autoref{fig.h=8}.

\subsection{Intersections of $h$-fragments}\label{s.intersections}
Given a graph~$\HH$, a subgraph $\Delta\subset\HH$ is called
\emph{perfect} if each vertex $u\in\HH\sminus\Delta$ is adjacent to at most
one vertex $v\in\Delta$. The \emph{perfect complement} of a perfect subgraph is
the set (see \autoref{S.6} or \autoref{S.8} for examples)
\[*
\HH\rminus\Delta:=\bigl\{u\in\HH\sminus\Delta\bigm|
 \text{$u$ is not adjacent to any $v\in\Delta$}\bigr\}.
\]

\proposition\label{prop.proper}
Let $\HH_1,\HH_2$ be two $h$-fragments in an admissible
graph~$\Gamma$. Then\rom:
\roster
\item
the intersection $\Delta:=\HH_1\cap\HH_2$ is a perfect subgraph of
each~$\HH_i$, $i=1,2$\rom;
\item
the relation $(u_1,u_2)\mapsto(u_1\cdot u_2=1)$ on
$(\HH_1\sminus\Delta)\times(\HH_2\sminus\Delta)$
is a bijection
$\att\:\HH_1\rminus\Delta\simeq\HH_2\rminus\Delta$, henceforth called the
\emph{attachment bijection}.
\endroster
\endproposition

\proof
Since each line intersects each hyperplane at one point (or because
$u\cdot h=1$ and $h=\sum_{v\in\HH}v$, \cf. \autoref{rem.admissible}), we have
\[
\text{each vertex $u\in\Gamma\sminus\HH$ is adjacent to exactly one
 vertex $v\in\HH$}.
\label{eq.one.vertex}
\]
If $\Delta$ were not perfect in, say, $\HH_1$, there would be a vertex
$u\in\HH_1\sminus\Delta=\HH_1\sminus\HH_2$ adjacent to at least two
vertices of $\Delta\subset\HH_2$, contradicting~\eqref{eq.one.vertex}.

Similarly, each vertex of $\HH_i\sminus\Delta$ must be adjacent to exactly
one vertex of~$\HH_j$, $j\ne i$. The vertices in
$(\HH_i\sminus\Delta)\sminus(\HH_i\rminus\Delta)$ are already adjacent
to those in $\Delta\subset\HH_j$, leaving a bijection between the perfect
complement.
\endproof

Our proof of \autoref{th.main} consists in listing all ``interesting'' graphs
by adding one $h$-fragment at a time. Thus, we start with an
admissible graph~$\Gamma_0$ and an $h$-fragment~$\HH$ and construct a new
graph $\Gamma:=\Gamma_0\cup\HH$, which must also be admissible.
For the new lattice $\Fano_{2d}(\Gamma)$, we need to choose and fix
\roster*
\item
the
intersection $\Delta:=\Gamma_0\cap\HH$ (or, rather, an isomorphism between
a pair of induced subgraphs $\Delta_0\subset\Gamma_0$ and
$\Delta_\HH\subset\HH$) and
\item
for each vertex $u\in\HH\sminus\Delta_\HH$, the set
$\Gamma_0(u):=\{v\in\Gamma_0\,|\,v\cdot u=1\}$.
\endroster

Since we need an admissible graph, these choices are quite limited, and it is
this fact that makes the computation feasible. First of all, till the very
last moment (\viz. saturating subgeometric graphs to geometric) we can
confine ourselves to the graphs~$\Gamma_0$ that are unions of $h$-fragments;
then, we can apply \autoref{prop.proper} to each pair $\HH,\HH_0$,
where $\HH_0\subset\Gamma_0$ is an $h$-fragment. Besides, we have the
following ``global'' analogue of \autoref{prop.proper}, whose proof is
straightforward, \cf.~\eqref{eq.one.vertex}.

\proposition\label{prop.global}
Let $\HH$ be an $h$-fragment and assume that the graph
$\Gamma:=\Gamma_0\cup\HH$ is admissible. Denote
$\Delta:=\Gamma_0\cap\HH$. Then\rom:
\roster
\item
$\Delta$ is a perfect subgraph of~$\Gamma_0$ \rom(but not necessarily
of~$\HH$\rom)\rom;
\item
the relation $(u_1,u_2)\mapsto(u_1\cdot u_2=1)$ on
$(\Gamma_0\sminus\Delta)\times(\HH\sminus\Delta)$
is a function
$\att\:\Gamma_0\rminus\Delta\to\HH\sminus\Delta$, henceforth called the
\emph{attachment map}.
\done
\endroster
\endproposition

\subsection{Idea of the proof of \autoref{th.main}}\label{s.idea}
As already explained, for degrees $2d\ge6$ we merely list all $2$-geometric
graphs~$\Gamma$ such that $\Fano_{2d}(\Gamma)$ is rationally generated by a
union of $h$-fragment. We use a modified version of the algorithm
of~\cite{degt:lines,degt.Rams:octics}, adding at each step a whole
$h$-fragment rather than a single vertex. We also use a number of enhancements
speeding up the algorithm; the  noteworthy ones are listed below.

\subsec\label{ss.limits}
At each step, the data
\[*
\Gamma_0\supset\Delta_0\simeq\Delta_\HH\subset\HH\quad\text{and}\quad
 \bigl\{\Gamma_0(u)\bigm|u\in\HH\sminus\Delta_\HH\bigr\}
\]
controlling the output (see \autoref{s.intersections}) are limited using
Propositions~\ref{prop.proper} and~\ref{prop.global}. In some cases we use
additional restrictions that are degree/graph specific,
\cf., for example, Lemmata~\ref{lem.h=6}
and~\ref{lem.h=8} below.

\subsec\label{ss.sat}
Each intermediate graph~$\Gamma$ obtained is replaced with its \emph{$h$-part}
\[*
\hp\Gamma:=\bigcup\HH,\quad
 \text{$\HH\subset\sat\Gamma$ is an $h$-fragment}
\]
generating the same lattice over~$\Q$; this makes \autoref{prop.proper}
more efficient.

\subsec\label{ss.iso}
After each step (and upon the normalization as in \autoref{ss.sat}), multiple
copies of the same graph are removed using the \texttt{digrap} package.

\subsec\label{ss.rank}
At each step we insist that the rank $\rank\Fano_{2d}(\Gamma)$ must increase,
as otherwise the new graph would be found among the saturations of the old
one.

\subsec\label{ss.rank.min}
For configurations of large rank $\rank\Fano_{2d}(\Gamma)\ge r_{\min}$ (a
threshold depending on the degree, typically $r_{\min}=18$), we switch to
adding one vertex at a time, literally as
in~\cite{degt:lines,degt.Rams:octics}.
In this case, \autoref{ss.sat} no longer applies.

\subsec\label{ss.other}
We also use a few tricks to control the order in which $h$-fragments are
added.
For example, typically $h$-fragments are added in the order
$\HH_1,\HH_2,\ldots$ in which they appear on the lists.
A few other degree specific details are discussed in the relevant sections
below.

\subsec\label{ss.config}
At the end, for further analysis, we compute all geometric graphs
containing at least a certain number $\#_{\min}$ of $h$-fragments. This
threshold
\[*
\minitab\quad
h^2=2d\:\    & 6&8&10&12&\ge14 \\
\#_{\min}\:\ &16&9& 3& 2& 1
\endminitab
\]
is chosen to keep the computation feasible and avoid large output.

\medskip
With this said, we consider \autoref{th.main} proved for all degrees $2d\ge6$.
\qed

\convention\label{conv.config}
When stating our results, depending on the output obtained, we
distinguish between \emph{configurations} (of lines) and
\emph{$h$-configurations}.

A \emph{configuration} is merely a geometric graph; it is called
\emph{maximal} if it is not a proper induced subgraph of another
configuration.

An $h$-configuration is a subgeometric graph~$\Xi$ with the following
properties:
\roster*
\item
$\Xi$ is a union (not necessarily disjoint) of $h$-fragments, and
\item
there is a geometric graph $\Gamma\supset\Xi$
such that $\Xi=\hp\Gamma$.
\endroster
It is \emph{maximal} if it is not a proper induced subgraph of another
$h$-configurations.

For example, the algorithm discussed in this section is mostly designed to
produce $h$-configurations. The counts stated in \autoref{th.main} are
those of $h$-configurations, not configurations.
\endconvention

\section{Quartics ($h\sp2=4$)}\label{S.4}

There is a unique $h$-fragment of degree~$4$, which is
\roster
\item\label{h=4,g=1}
  the complete graph $K(4)$.
\endroster
Obviously, $\rank\Z\HH_1=4$, $\girth\HH_1=3$, and $\Aut\HH_1=\SG4$.
Our principal result in this section is the following theorem.

\theorem\label{th.h=4}
The number of $h$-fragments on a smooth quartic $X\subset\Cp3$ is
\roster*
\item
$72$, if $X$ is Schur's quartic $X_{64}$, see \cite{Schur:quartics,DIS}\rom;
\item
$60$, if $X$ has $60$ lines, \ie, is
$X_{60}'$, $X_{60}''$, or $\bar X_{60}''$ in~\cite{DIS}\rom;
\item
less than $60$ otherwise.
\endroster
\endtheorem

\remark\label{rem.h=4}
We know a single example of a quartic with $50$ $h$-fragments, \viz. $Y_{56}$
in~\cite{DIS}; in all other known cases, the number of $h$-fragments is at
most $48$. Observe that $Y_{56}$ is also an example of a \emph{real} quartic
with $50$ \emph{real} $h$-sections.
\endremark

\conjecture\label{conj.h=4}
With the exception of the four surfaces mentioned in \autoref{th.h=4} and
\autoref{rem.h=4}, a smooth quartic $X\subset\Cp3$ has at most $48$
$h$-fragments.
\endconjecture


\proof[Proof of \autoref{th.h=4}]
All quartics~$X$ with $\ls|\Fn X|>48$ are known,
see~\cite{degt.Rams:quartics}; we check the $26$ graphs one-by-one and find $72$, $60$
(as in \autoref{th.h=4}), $50$ (as in \autoref{rem.h=4}), or at most $48$
$h$-fragments. Therefore, henceforth we assume that $\ls|\Fn X|\le48$.

As is well known (see, \eg,~\cite{DIS}), the star of a line~$\ell$ in
the Fano graph of a smooth quartic surface $X\subset\Cp3$ is of the form
$p\tA_2\oplus q\bA_1$, where the values of $(p,q)$ are as in the following
table:
\[*
\minitab\quad
p={}&   0& 1&2&3&4&5&6\\
q\le{}&12&10&9&7&6&3&2
\endminitab
\]
It follows that $\mult\ell\le6$ for each $\ell\in\Fn X$ and
we can easily bound $\Hnum(X)$ by~$72$: since
each $h$-fragment consists of $4$ lines, we have
\[
\Hnum(X)\le\frac{\max(\mult\ell)}{4}\ls|\Fn X|\le\frac64\cdot48=72.
\label{eq.h=4}
\]
To improve this bound, we follow~\cite{DIS} and compute all configurations
with a bouquet of six $h$-fragments, finding no new
surface~$X$ with $\Hnum(X)>48$. Hence, we can assume $\max(\mult\ell)\le5$
in~\eqref{eq.h=4},
arriving at $\Hnum(X)\le60$.

Finally, under the assumptions, for $\Hnum(X)=60$ we must have $\ls|\Fn X|=48$ and
$\mult\ell=5$, hence $\val\ell\ge15$, for \emph{each} line $\ell$.
Then, mimicking the proof of Segre's theorem~\cite{Segre}
cited in the first paragraph of
\autoref{S.intro}, we arrive at $\ls|\Fn X|\ge52$.
\endproof

\section{Sextics ($h\sp2=6$)}\label{S.6}

A smooth $K3$-sextic $X\subset\Cp4$ is a regular complete intersection of a
quadric and a cubic, $X=Q_2\cap Q_3$.
There are two $h$-fragments of degree~$6$, \viz.
\roster
\item\label{h=6,g=1}
  $\GRAPH\AA[2](1;2;3) \(1x2)\(1x3)\(2x3)$ (the $3$-prism) and
\item\label{h=6,g=2}
  the complete bipartite graph $K(3,3)$ (the utility graph),
\endroster
see also \autoref{fig.h=6} and \autoref{tab.h=6}.
Geometrically (see~\autoref{S.intro}), cutting $X=Q_2\cap Q_3$ by a
hyperplane~$H$, the two split sections differ by whether $H\cap Q_2$ is a smooth
quadric $\Cp1\times\Cp1$ or a pair of planes $\Cp2\cup\Cp2$.


\figure
\[*
\llap{\iref{h=6,g=1}}
\PIC(60,45)(0,0)
\(30,45){\line(2,-3){30}}
\(30,45){\line(-2,-3){30}}
\(30,45){\line(0,-1){15}}
\(30,30){\line(2,-3){15}}
\(30,30){\line(-2,-3){15}}
\(0,0){\line(1,0){60}}
\(0,0){\line(2,1){15}}
\(60,0){\line(-2,1){15}}
\(15,7.5){\line(1,0){30}}
\(30,45)\b\no(25,45)1
\(30,30)\b\no(30,20)4
\(0,0)\b\no(-6,-2)2
\(60,0)\b\no(66,-2)3
\(15,7.5)\b\no(22,9)5
\(45,7.5)\b\no(38,9)6
\endPIC
\qquad\qquad
\llap{\iref{h=6,g=2}\ \ }
\PIC(60,30)(0,0)
\(0,30){\line(0,-1){30}}
\(0,30){\line(1,-1){30}}
\(0,30){\line(2,-1){60}}
\(60,30){\line(0,-1){30}}
\(60,30){\line(-1,-1){30}}
\(60,30){\line(-2,-1){60}}
\(30,30){\line(0,-1){30}}
\(30,30){\line(1,-1){30}}
\(30,30){\line(-1,-1){30}}
\(0,0)\b\no(0,-8)4
\(30,0)\b\no(30,-8)5
\(60,0)\b\no(60,-8)6
\(0,30)\b\no(0,34.5)1
\(30,30)\b\no(30,34.5)2
\(60,30)\b\no(60,34.5)3
\endPIC
\]
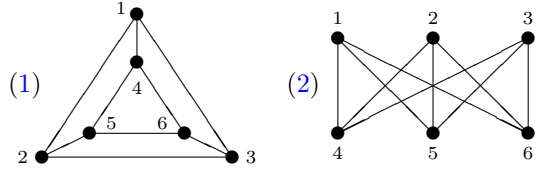
\caption{The $h$-fragments of degree $6$}\label{fig.h=6}
\endfigure

\table
\caption{$h$-configurations of degree~$6$ \noaux{(see \autoref{conv.tab1})}}\label{tab.h=6}
\vspace{-2\bigskipamount}
\[*
\let\1\tabconfig
\minitab{\kern6pt}
\HHH&(r,g,s)&\omit\hss$h$-fragment counts\hss
                       &&\multispan5\minitabskip\hss large counts\hss\\\trule
\iref{h=6,g=1}&(6,3,12) &\clist{0..6, 8..10, 12, 15, 16, 18, 20, 25, 36}
                       &&16&\.&36&25&36\\
\iref{h=6,g=2}&(6,4,72) &\clist{0..24, 26, 27, 29, 30, 36, 40}
                       &&20&36&12&24&40\\
\counts{total counts}   &\clist{0..34, 36, 48, 49, 76}
                       &&36&36&48&49&76\\
\counts{configuration}& &&34&36&36&\1{\Psi_{38}}&\1{\Psi_{42}}
\endminitab
\]
\endtable

\convention\label{conv.tab1}
Listed in Tables~\ref{tab.h=6}, \ref{tab.h=8}, \ref{tab.h=10},
are the values taken by the number of $h$-fragments,
both total and itemized by the type of the graphs, and, in the last column,
the itemized counts for a few largest $h$-configurations. In the last row,
for each of these few largest $h$-configurations, we show the corresponding (minimal)
configuration, either by its ``name'' in
\cite{degt:lines,DR:octic.graphs,degt.Rams:octics} or just by the number of
lines.
\endconvention

As in the case of
quartics (see \autoref{th.h=4}), the two largest configurations
$\Psi_{38}$, $\Psi_{42}$ maximize also the number of lines,
see~\cite{degt:lines}; they are triangular (or $\tA_2$-, see
\autoref{s.fragments}).
The other three are quadrangular (or $\tA_3$-, see
\autoref{s.fragments});
they can be identified by the triples
$(\ls|\Gamma|,\rank\Z\Gamma,\ls|\Aut\Gamma|)$,
which are as
follows:
\[*
(16+20)\mapsto(34,18,32),\quad
(0+36)\mapsto(36,19,24),\quad
(36+12)\mapsto(36,18,1296).
\]
Probably, they will be given ``names'' in~\cite{degt.Rams:sextics}.

Our original strategy for proving the bound $\Hnum(X)\le76$ was bounding the
multiplicities $\mult\ell$ first and then arguing as in \autoref{S.4}.
We found that
this requires almost as much computation as compiling the full
list of $h$-configurations. Nevertheless, we discuss a few simple facts
concerning bouquets of $h$-fragments of degree~$6$.

Ultimately,
denoting by~$c_i$ the number of copies of $\HH_i$, we have
\[*
\minitab\quad
c_2=  &0&1&2&3&4&5&6&7&8&9&10\\
c_1\le&5&6&6&5&5&6&6&4&5&0&0\\
\endminitab
\]
If $c_2=8$, then $c_1\in\{0, 2, 5\}$.
To describe the structure of a bouquet, recall (see
\cite[\S\,2.5]{degt:lines}) that the star of a vertex $\ell\in\Fn X$ is of the
form
\[
p\bA_2\oplus q\bA_1,\qquad
 p\le1,\quad
 q\le9,\quad
 q\ne8\ \text{if}\ p=1.
\label{eq.h=6}
\]
By \autoref{prop.proper}, two $h$-fragments intersects in a perfect
subset~$\Delta$. Up
to automorphism, proper perfect subsets are: in $\HH_1$,
\roster
\item\label{bq.1.1}
$\Delta=\{1\}$, with $\HH_1\rminus\Delta=\{5,6\}$;
\item\label{bq.1.2}
$\Delta=\{1,4\}$, with $\HH_1\rminus\Delta=\varnothing$;
\item\label{bq.1.3}
$\Delta=\{1,2,3\}$, with $\HH_1\rminus\Delta=\varnothing$,
\endroster
and in $\HH_2$,
\roster[\lastitem]
\item\label{bq.2.1}
$\Delta=\{1\}$, with $\HH_2\rminus\Delta=\{2,3\}$;
\item\label{bq.1.2}
$\Delta=\{1,4\}$, with $\HH_2\rminus\Delta=\varnothing$.
\endroster
Here and below, we index the vertices $v_k\in\HH_i$ as shown in the figures,
and we often abbreviate~$v_k$ to just $k$.

The following lemma can be used to speed up the computation, see
\autoref{ss.limits}.

\lemma\label{lem.h=6}
The intersection $\Delta:=\HH_1\cap\HH_1'$ of two copies of $\HH_1$ is either
empty or a common triangle.
\endlemma

\proof
If $\HH_1$ contains a line~$\ell$, it must also contain an $\bA_2$-subgraph
in the star of~$\ell$. According to~\eqref{eq.h=6}, this subgraph is unique.
\endproof

We arive at the following description of a bouquet. Let $S$ be the set of the
$\bA_1$-components of the star of~$\ell$; we have $\ls|S|\le9$,
see~\eqref{eq.h=6}.
Replacing each $h$-fragment $\HH\ni\ell$ with its \emph{germ} at~$\ell$,
\ie, the set of its vertices that are adjacent to~$\ell$, we get a set
$\fS_1\subset S$ (the germs of $\HH_1$-fragments without the common
$\bA_2$-type component) and a set~$\fS_2$ of $3$-element subsets of~$S$ (the
germs of $\HH_2$-fragments). The latter has the property that
\[*
\ls|s'\cap s''|\le1\quad\text{for any pair $s',s''\in\fS_2$, $s'\ne s''$}.
\]
The attachment bijections are less straightforward, and we omit their
description.

\section{Octics ($h\sp2=8$)}\label{S.8}

A typical smooth octic $K3$-surface $X\subset\Cp5$ is a \emph{triquadric}, \ie,
regular complete intersection of three quadric hypersurfaces. However, in the
space of all octics there is a divisor of \emph{special} octics which cannot
be represented in this way. According to \cite{Saint-Donat} (see also
\cite{degt.Rams:octics} for a homological restatement), a smooth octic~$X$ is
a triquadric if and only if the polarized lattice $\NS(X)\ni h$ is
$3$-admissible.

There are three $h$-fragments of degree~$8$:
\roster
\item\label{h=8,g=1}
  $\GRAPH\AA[2](1;2;3) \A[1](1,2,3) \A[1](1,2,3) $;
\item\label{h=8,g=2}
  $\GRAPH\AA[3](1;2;3;4) \(1x2)\(1x3)\(2x4)\(3x4)$ (the Wagner graph);
\item\label{h=8,g=3}
  $\GRAPH\AA[3](1;2;3;4) \(1x2)\(1x4)\(2x3)\(3x4)$ (the $3$-cube),
\endroster
see also \autoref{fig.h=8} and \autoref{tab.h=8} (where $\HH_1$ and
$\HH_2,\HH_3$ are treated separately due to \autoref{prop.h=8} below).
\figure
\[*
\llap{\iref{h=8,g=1}}
\PIC(60,45)(0,0)
\(30,45){\line(2,-3){30}}
\(30,45){\line(-2,-3){30}}
\(30,45){\line(0,-1){15}}
\(30,30){\line(2,-3){15}}
\(30,30){\line(-2,-3){15}}
\(0,0){\line(1,0){60}}
\(0,0){\line(2,1){37}}
\(60,0){\line(-2,1){37}}
\(30,45)\b
\(30,45)\c
\(30,30)\b
\(0,0)\b
\(0,0)\c
\(60,0)\b
\(60,0)\c
\(15,7.5)\b
\(45,7.5)\b
\(37.5,18.75)\b
\(22.5,18.75)\b
\endPIC
\qquad\qquad
\def\cube{
\(0,0)\b\no(-6,-2)0
\(45,0)\b\no(51,-2)1
\(0,45)\b\no(-6,43)4
\(45,45)\b\no(51,43)5
\(12,12)\b\no(15,5)2
\(12,33)\b\no(15,36)6
\(33,33)\b\no(30,36)7
\(33,12)\b\no(30,5)3
}
\llap{\iref{h=8,g=2}\ \ }
\PIC(45,45)(0,0)
\(0,45){\line(0,-1){45}}
\(0,45){\line(1,0){45}}
\(0,45){\line(1,-1){45}}
\(0,0){\line(1,0){45}}
\(0,0){\line(1,1){45}}
\(45,0){\line(0,1){45}}
\(12,33){\line(1,0){21}}
\(33,12){\line(-1,0){21}}
\cube
\endPIC
\qquad\qquad
\llap{\iref{h=8,g=3}\ \ }
\PIC(45,45)(0,0)
\(0,45){\line(0,-1){45}}
\(0,45){\line(1,0){45}}
\(0,45){\line(1,-1){12}}
\(0,0){\line(1,0){45}}
\(0,0){\line(1,1){12}}
\(45,0){\line(0,1){45}}
\(45,0){\line(-1,1){12}}
\(45,45){\line(-1,-1){12}}
\(12,33){\line(0,-1){21}}
\(12,33){\line(1,0){21}}
\(33,12){\line(0,1){21}}
\(33,12){\line(-1,0){21}}
\cube
\endPIC
\]
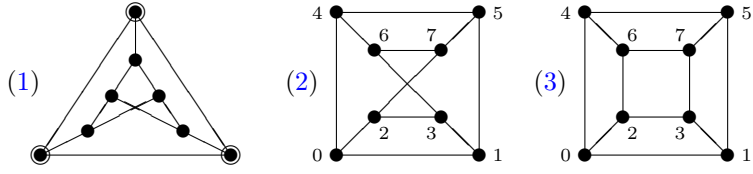
\caption{The $h$-fragments of degree $8$}\label{fig.h=8}
\endfigure

\table
\caption{$h$-configurations of degree~$8$ \noaux{(see \autoref{conv.tab1})}}\label{tab.h=8}
\vspace{-2\bigskipamount}
\[*
\let\1\tabconfig
\minitab{\kern6pt}
\HHH&(r,g,s)&\omit\hss$h$-fragment counts\hss
                       &&\multispan7\minitabskip\hss large counts\hss\\\trule
\iref{h=8,g=1}&(8,3,12) &\clist{0..10, 12, 14, 16, 18, 20, 36, 42, 72}
                       &&72\\
\iref{h=8,g=2}&(8,4,16) &\clist{0..16, 18..21, 24, 32, 48}
                       &&&24&16&48&32&48&\.\\
\iref{h=8,g=3}&(8,4,48) &\clist{0..13, 16, 20, 21, 32, 80}
                       &&&12&20& 8&32&16&80\\
\counts{total $\iref{h=8,g=2}+\iref{h=8,g=3}$}
                        &\clist{0..18, 20, 21, 23, 24, 26, 36, 56, 64, 80}
                       &&&36&36&56&64&64&80\\
\counts{configuration}&&&\1{\Psi_{33}}&28&28&\1{\Theta_{32}'}&\1{\Theta_{32}''}&\1{\Theta_{32}'''}&\1{\Theta_{32}^\mathrm{K}}
\endminitab
\]
\endtable

The maximal number~$72$ of $\HH_1$-type fragments is attained at the
configuration $\Psi_{33}$ maximizing the number of lines on \emph{special}
octics, see~\cite[Theorem~1.3]{degt:lines}. Unlike the cases of quartics and
sextics, the Kummer triquadric $\Theta_{32}^\mathrm{K}$ with the maximal
number of $h$-fragments is not line maximizing (and all its fragments are of the
same type).
One of the two line-maximizing configurations, \viz.
$\Theta_{36}'\supset\Theta_{34}'\supset\Theta_{32}'$,
has but $56$ $h$-fragments; the other one, $\Theta_{36}''$, has no $h$-fragments
at all.

\proposition\label{prop.h=8}
The fragment~$\HH_1$ appears in special octics only, whereas $\HH_2$,
$\HH_3$ appear in triquadrics only.
\endproposition

\proof
The fragment $\HH_1$ has a
$3$-isotropic vector, \eg, the sum of the three vertices
(circled in \autoref{fig.h=8}) of the outer
triangle.
For $\HH=\HH_2$ or~$\HH_3$, we try to append a $3$-isotropic vector~$p$, \ie,
consider the lattice
\[*
\Cal{N}:=(\Z\HH+\Z h+\Z p)/\!\ker,\qquad
 p^2=0,\quad p\cdot h=3.
\]
This notation is ambiguous as we do not know the intersections of~$p$ with the
vertices of~$\HH$.
By \cite[Lemma 2.29]{degt.Rams:octics}, and since $h=\sum_{v\in\HH}v$, there
is a $3$-element subset $o\subset\HH$ such that $p\cdot v=1$ if $v\in o$ and
$p\cdot v=0$ otherwise. We try all such subsets and find that the polarized
lattice $\Cal{N}\ni h$ is not $2$-admissible.

There are five $\Aut\HH_2$-orbits of
$3$-element subsets $o\subset\HH_2$:
\roster*
\item
$o=\{ 0, 1, 2 \}$ or $\{ 0, 2, 4 \}$:
the lattice is not hyperbolic, $\Gs_+\Cal{N}=2$;
\item
$o=\{ 0, 1, 6 \}$: there is an exceptional divisor $v_2+v_3+v_7-p$;
\item
$o=\{ 0, 2, 5 \}$: there is an exceptional divisor $v_3+v_4+v_6-p$;
\item
$o=\{ 0, 3, 5 \}$: there is an exceptional divisor $v_4+v_6+v_7-p$,
\endroster
and three $\Aut\HH_3$-orbits of
$3$-element subsets $o\subset\HH_3$:
\roster*
\item
$o=\{ 0, 1, 2 \}$:
the lattice is not hyperbolic, $\Gs_+\Cal{N}=2$;
\item
$o=\{ 0, 1, 6 \}$: there is an exceptional divisor $v_2+v_3+v_7-p$;
\item
$o=\{ 0, 3, 5 \}$: there is an exceptional divisor $v_4+v_6+v_7-p$.
\endroster
This concludes the proof.
\endproof

\remark
In the proof, the exceptional divisor is of the form
$-p+\sum_{v\in\Xi}v$, where $\Xi\subset\HH$ is an induced subgraph
isomorphic to $\bA_3$ and disjoint from~$o$.
Such a divisor exists for any $3$-element set $o\subset\HH$;
it can be used instead of $\Gs_+\Cal{N}$.
\endremark

\remark\label{rem.h=8}
According to \cite[\S\,2.5]{degt:lines}, one has $\val\ell\le7$ for each
vertex~$\ell$ of the Fano graph $\Fn X$ of a smooth triquadric~$X$, and
configurations with $7$-valent vertices do appear on our list of
configurations, \eg, $\Theta'_{36}\supset\Theta'_{34}$.
Nevertheless, experimentally we observe that
$\val\ell\le6$ \emph{in the $h$-part $\hp(\Fn X)$}.
Alternatively, this fact could be proved by considering
the few dozens of minimal bouquets whose $h$-fragments together
make use
of all seven vertices adjacent to~$\ell$, \cf. the discussion of bouquets
below.
We omit this routine search for exceptional divisors.
\endremark

\problem\label{prob.h=8}
Geometrically, $\HH_2$ and $\HH_3$ represent split degree~$8$ curves serving as base
loci of certain nets of quadrics in $\Cp4$. It appears that $\HH_3$ is
generic whereas $\HH_2$ is degenerate, even though they have the same
codimension, but we do not know a precise explanation of this
phenomenon.
\endproblem

The number of isomorphism classes of bouquets ($837$)
is almost the same as that of configurations
($860$). It appears that only bouquets of $\HH_3$-fragments have a nice
description. The proper perfect subsets of $\HH_3$ are, up to automorphism,
\roster
\item\label{bq.0}
$\Delta=\{0\}$, with $\HH_3\rminus\Delta=\{3,5,6,7\}$;
\item\label{bq.edge}
$\Delta=\{0,2\}$, with $\HH_3\rminus\Delta=\{5,7\}$;
\item\label{bq.0-7}
$\Delta=\{0,7\}$, with $\HH_3\rminus\Delta=\varnothing$;
\item\label{bq.face}
$\Delta=\{0,1,2,3\}$, with $\HH_3\rminus\Delta=\varnothing$,
\endroster
and we can assume that $v_0$ in each fragment is the central line of the
bouquet.

The following lemma speeds up the computation, see \autoref{ss.limits}.

\lemma\label{lem.h=8}
There are the following extra restrictions on $\Delta:=\HH_3\cap\HH_3'$\rom:
\roster*
\item
$\ls|\Delta|\ne1$, \ie, case~\iref{bq.0} is not realizable\rom;
\item
in case~\iref{bq.edge}, the attachment bijection
$\HH_3\rminus\Delta\to\HH_3'\rminus\Delta$ is
$(5,7)\leftrightarrow(7',5')$.
\endroster
\endlemma

\proof
In case~\iref{bq.0}, up to automorphisms of the two graphs,
the attachment bijection is either
$\att\:(3,5,6,7)\leftrightarrow(3',5',6',7')$ or
$(3,5,6,7)\leftrightarrow(3',5',7',6')$; each
results in a number of exceptional
divisors, \eg,
\[*
e=-h+v_4+v_5+v_6+v_7+v'_4+v'_5+v'_6+v'_7.
\]
In case~\iref{bq.edge}, the bijection
$\att\:(5,7)\leftrightarrow(5',7')$ results in the exceptional divisor
\[*
e=-h+2v_2+v_3+v_6+v_7+v'_3+v'_6+v'_7,
\]
leaving
$\att\:(5,7)\leftrightarrow(7',5')$ as the only choice.
\endproof

Hence, we have the following simple description of a
bouquet of $\HH_3$-fragments at a line~$\ell$. Let $S$ be the set of lines
adjacent to~$\ell$; by \autoref{rem.h=8}, we can assume
$\ls|S|\le6$. Replacing each $\HH_3$-fragment with its germ at~$\ell$, we
obtain a set $\fS_3$ of $3$-element subsets $s\subset S$, and the rest of the
graph is recovered from~$\fS_3$ by \autoref{lem.h=8}.

There are no immediate restrictions on~$\fS_3$, as the set of all
$3$-element subsets gives rise to a geometric graph. Hence, all other graphs
are subgeometric, but not all of them are saturated. Altogether, there are
$14$ isomorphism classes of saturated pairs $(\text{bouquet},\ell)$;
their numbers of $\HH_3$-fragments are
\[*
 1, 2, 2, 2, 3, 3, 4, 4, 4, 5, 5, 7, 8, 20.
\]

\remark\label{rem.Theta32}
The maximal bouquet with $20$ copies of $\HH_3$ is the graph
$\Theta^\mathrm{K}_{32}$ that maximizes the number of $h$-fragments. In this
configuration, all $32$ bouquets are isomorphic, each consisting of $20$ fragments.
\endremark

\remark
As a technical detail, in the computation for triquadrics the graphs are
processed in the order $\HH_3,\HH_2$, see \autoref{ss.other}.
\endremark

\section{The polarization $h\sp2=10$}\label{S.10}

There are six $h$-fragments of degree~$10$:
\roster
\item\label{h=10,g=1}
  $\GRAPH\AA[2](1;2;3) \AA[2](1;2;3) \A[1](1,2,3) $;
\item\label{h=10,g=2}
  $\GRAPH\AA[3](1;2;3;4) \A[2](1,2;3,4) \(1x3)\(2x4)$;
\item\label{h=10,g=3}
  $\GRAPH\AA[3](1;2;3;4) \A[1](1,2,4) \A[1](2,3,4) \(1x3)$;
\item\label{h=10,g=4}
  $\GRAPH\AA[3](1;2;3;4) \A[2](1,3;2,4) \(1x2)\(3x4)$;
\item\label{h=10,g=5}
  $\GRAPH\AA[3](1;2;3;4) \A[2](1,4;2,3) \(1x2)\(3x4)$;
\item\label{h=10,g=6}
  $\GRAPH\AA[4](1;2;3;4;5) \(1x3)\(1x4)\(2x4)\(2x5)\(3x5)$ (Petersen graph),
\endroster
see also \autoref{fig.h=10} (where only two graphs are shown,
merely to illustrate how weird
these graphs can be; though, \cf. also \autoref{rem.h=10} below)
and \autoref{tab.h=10}.

Among the three
line-maximizing configurations only one, \viz. $\Phi_{30}'''$
(see~\cite{degt:lines}) has many $h$-fragments, \viz. $14$ copies of
the Petersen graph $\HH_6$. The other two have no $h$-fragments at all.

\figure
\[*
\llap{\iref{h=10,g=1}}
\PIC(60,45)(0,0)
\(30,45){\line(2,-3){30}}
\(30,45){\line(-2,-3){30}}
\(30,45){\line(0,-1){30}}
\(30,30){\line(-2,-3){7.5}}
\(0,0){\line(1,0){60}}
\(0,0){\line(2,1){30}}
\(60,0){\line(-2,1){30}}
\(15,7.5){\line(1,0){15}}
\(45,7.5){\line(-2,3){7.5}}
\(22.5,18.75){\line(1,0){15}}
\(22.5,18.75){\line(2,-3){7.5}}
\(37.5,18.75){\line(-2,-3){7.5}}
\(30,45)\b
\(30,45)\c
\(30,30)\b
\(0,0)\b
\(0,0)\c
\(60,0)\b
\(60,0)\c
\(15,7.5)\b
\(45,7.5)\b
\(30,7.5)\b
\(30,7.5)\c
\(37.5,18.75)\b
\(37.5,18.75)\c
\(22.5,18.75)\b
\(22.5,18.75)\c
\(30,15)\b
\endPIC
\ \cdots\quad\qquad
\llap{\iref{h=10,g=6}\ \ }
\PIC(60,53)(0,0)
\(0,33){\line(3,2){30}}
\(0,33){\line(1,-3){11}}
\(0,33){\line(4,-1){15}}
\(60,33){\line(-3,2){30}}
\(60,33){\line(-1,-3){11}}
\(60,33){\line(-4,-1){15}}
\(11,0){\line(1,0){38}}
\(11,0){\line(1,1){10}}
\(21,10){\line(6,5){24}}
\(21,10){\line(1,3){9}}
\(49,0){\line(-1,1){10}}
\(39,10){\line(-6,5){24}}
\(39,10){\line(-1,3){9}}
\(30,53){\line(0,-1){16}}
\(15,29){\line(1,0){30}}
\(15,29)\b
\(45,29)\b
\(0,33)\b
\(60,33)\b
\(30,53)\b
\(30,37)\b
\(11,0)\b
\(21,10)\b
\(39,10)\b
\(49,0)\b
\endPIC
\]
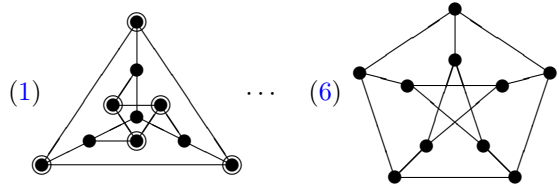
\caption{Some $h$-fragments of degree $10$}\label{fig.h=10}
\endfigure

\table
\caption{$h$-configurations of degree~$10$ \noaux{(see \autoref{conv.tab1})}}\label{tab.h=10}
\vspace{-2\bigskipamount}
\[*
\let\1\tabconfig
\minitab{\kern6pt}
\HHH&(r,g,s)&\omit\hss$h$-fragment counts\hss
                       &&\multispan9\minitabskip\hss large counts\hss\\\trule
\iref{h=10,g=1}&(9,3,12) &\clist{0, 1, 3, 6, 15 }
                       &&\.&\.&\.&\.&\.&\.&\.&15&\.\\
\iref{h=10,g=2}&(10,4,4) &\clist{0..6}
                       &&\.&\.&\.&\.&\.&\.&\.&\.&\.\\
\iref{h=10,g=3}&(10,4,8) &\clist{0..6, 8 }
                       && 6&\.&\.&\.&\.&\.&\.&\.&\.\\
\iref{h=10,g=4}&(10,4,20) &\clist{0..4, 6, 8, 12 }
                       && 6&12& 8& 6& 6&\.&\.&\.&\.\\
\iref{h=10,g=5}&(10,4,20) &\clist{0..4, 6, 12}
                       &&\.&\.& 4& 6& 6&12&\.&\.&\.\\
\iref{h=10,g=6}&(10,5,120) &\clist{0..8, 11, 14, 16}
                       &&\.&\.&\.&\.&\.& 2&14&\.&16\\
\counts{total counts}     &\clist{0..12, 14..16}
                       &&12&12&12&12&12&14&14&15&16\\
\counts{configuration}&&&26&20&20&20&20&20&\1{\Phi_{30}'''}&22&28
\endminitab
\]
\endtable

\remark\label{rem.h=10}
This is the first degree where the strata of surfaces with a fixed
$h$-fragment $\HH_i$ may have different codimension $r:=\rank\HH_i$.
In a sense, $\HH_1$ is ``less degenerate'' than the others.
Geometrically, this hyperplane section contains two parabolic fibers of the
same elliptic pencil, \viz. the two triangles with circled vertices in
\autoref{fig.h=10}, left.
It is worth mentioning that this special graph~$\HH_1$ never appears together
with any other $h$-fragment in the same configuration.
\endremark

\remark\label{rem.pentagon}
As a technical detail (see \autoref{ss.other}),
when the first copy of $\HH_6$ is added to a graph~$\Gamma_0$ containing
$\HH_1,\ldots,\HH_5$ only,
\emph{and if $\Delta:=\HH_6\cap\Gamma_0$ is to be empty},
we proceed in two steps: first, a pentagon $\tA_4$ is added, and then the
rest of~$\HH_6$.
For some graphs~$\Gamma_0$
this trick speeds up the computation.
\endremark

\section{The polarization $h\sp2=12$}\label{S.12}

There are nine $h$ fragments of degree~$12$:
\roster
\item\label{h=12,g=1}
  $\GRAPH\AA[2](1;2;3) \AA[2](1;2;3) \AA[2](1;2;3) $;
\item\label{h=12,g=2}
  $\GRAPH\AA[3](1;2;3;4) \AA[3](1;2;3;4) \(1x2)\(3x4)$;
\item\label{h=12,g=3}
  $\GRAPH\AA[3](1;2;3;4) \AA[3](1;2;4;3) \(1x2)\(3x4)$;
\item\label{h=12,g=4}
  $\GRAPH\AA[3](1;2;3;4) \AA[3](1;2;3;4) \(1x3)\(2x4)$;
\item\label{h=12,g=5}
  $\GRAPH\AA[3](1;2;3;4) \A[2](1,2;3,4) \A[2](1,3;2,4) $;
\item\label{h=12,g=6}
  $\GRAPH\AA[3](1;2;3;4) \A[2](1,2;3,4) \A[2](1,4;2,3) $;
\item\label{h=12,g=7}
  $\GRAPH\AA[3](1;2;3;4) \A[2](1,3;2,4) \A[2](1,3;2,4) $;
\item\label{h=12,g=8}
  $\GRAPH\AA[4](1;2;3;4;5) \A[1](1,4,5) \A[1](2,3,5) \(1x3)\(2x4)$;
\item\label{h=12,g=9}
  $\GRAPH\AA[4](1;2;3;4;5) \A[1](1,2,5) \A[1](3,4,5) \(1x3)\(2x4)$,
\endroster
see also \autoref{tab.h=12}.
We conclude that $\Hnum(X)\le90$ is the sharp bound.

\table
\caption{$h$-configurations of degree~$12$ \noaux{(see \autoref{conv.tab4})}}\label{tab.h=12}
\vspace{-2\bigskipamount}
\[*
\def\1{}
\minitab{\kern4.2pt}
\HHH&(r,g,s)&\multispan{28}\ \hss$h$-fragment counts\hss\\\trule
\iref{h=12,g=1}&(10,3,36) &20&4 &1 &\.&\.&\.&\.&\.&\.&\.&\.&\.&\.&\.&\.&\.&\.&\.&\.&\.&\.&\.&\.&\.&\.&\.&\.&\.\\
\iref{h=12,g=2}&(10,4,24) &\.&\.&\.&1 &\.&\.&\.&\.&\.&\.&\.&\.&\.&\.&\.&\.&\.&\.&\.&\.&\.&\.&\.&\.&\.&\.&\.&\.\\
\iref{h=12,g=3}&(10,4,24) &\.&\.&\.&\.&1 &\.&\.&\.&\.&\.&\.&\.&\.&\.&\.&\.&\.&\.&\.&\.&\.&\.&\.&\.&\.&\.&\.&\.\\
\iref{h=12,g=4}&(11,4,16) &\.&\.&\.&\.&\.&3 &1 &\.&\.&\.&\.&\.&\.&\.&\.&\.&\.&\.&\.&\.&\.&\.&\.&\.&\.&\.&\.&\.\\
\iref{h=12,g=5}&(12,4,4)  &\.&\.&\.&\.&\.&\.&\.&16&4 &3 &2 &1 &1 &\.&\.&\.&\.&\.&\.&\.&\.&\.&\.&\.&\.&\.&\.&\.\\
\iref{h=12,g=6}&(12,4,8)  &\.&\.&\.&\.&\.&\.&\.&\.&\.&\.&\.&1 &\.&6 &1 &\.&\.&\.&\.&\.&\.&\.&\.&\.&\.&\.&\.&\.\\
\iref{h=12,g=7}&(12,4,48) &\.&\.&\.&\.&\.&\.&\.&\.&\.&\.&\.&\.&\.&\.&\.&16&2 &1 &\.&\.&\.&\.&\.&\.&\.&\.&\.&\.\\
\iref{h=12,g=8}&(12,5,16) &\.&\.&\.&\.&\.&\.&\.&\.&\.&\.&\.&\.&\.&\.&\.&\.&\.&\.&90&10&6 &4 &3 &2 &1 &\.&\.&\.\\
\iref{h=12,g=9}&(12,5,18) &\.&\.&\.&\.&\.&\.&\.&\.&\.&\.&\.&\.&\.&\.&\.&\.&\.&\.&\.&\.&\.&\.&\.&\.&\.&3 &2 &1 \\
\counts{$h$-configs}      &\1&\1&\1&\1&\1&\1&\1&\1&2 &\1&6 &\1&\1&\1&\1&\1&2 &\1&\1&\1&\1&3 &3 &5 &\1&\1&2 &\1\\
\endminitab
\]
\endtable

\convention\label{conv.tab4}
In Tables~\ref{tab.h=12}--\ref{tab.h=18}, we show the itemized $h$-fragment
conunts
for each $h$-configuration, listed one-by-one (except \autoref{tab.h=12},
where $h$-configurations with the same conuts are combined; their number, if
greater than~$1$, is given in the last row).
In Tables~\ref{tab.h=16} and~\ref{tab.h=18}, the last row is the number of
configurations.
\endconvention

The four largest $h$-configurations are as follows:
\roster*
\item
$90\times\HH_8$ is the only line-maximizing graph $\Phi''_{36}$,
see~\cite{degt:lines}; its automorphism group is transitive, each bouquet
consisting of $30$ $h$-fragments;
\item
$20\times\HH_1$ is a certain maximal triangular configuration~$\Gamma$ with $21$ lines,
not transitive;
one has $\rank\Z\Gamma=14$ and $\ls|\Aut\Gamma|=4320$.
\item
$16\times\HH_7$ is a certain maximal quadrangular configuration~$\Gamma$ with $24$
lines;
it is transitive, each bouquet consisting of eight $h$-fragments, and
one has $\rank\Z\Gamma=14$ and $\ls|\Aut\Gamma|=768$;
\item
$16\times\HH_5$ is a certain quadrangular configuration~$\Gamma$ with $24$
lines, neither transitive nor maximal;
one has $\rank\Z\Gamma=17$ and $\ls|\Aut\Gamma|=64$.
\endroster

\remark
As a technical detail, we act as in \autoref{rem.pentagon}
when the first copy of $\HH_9$ is added to a graph~$\Gamma_0$
and if $\Delta:=\HH_6\cap\Gamma_0$ is to be empty.
\endremark

\section{Higher degrees}\label{S.higher}

Starting from degree $h^2=14$, we find all \emph{configurations}
containing at least one $h$-fragment, that is,
starting
the computation with an
$h$-fragment~$\HH$,
we immediately switch to adding one vertex at
a time, see \autoref{ss.rank.min}.

If $2d\ge14$, all $h$-configurations are geometric.

\subsection{The case $h\sp2=14$}\label{s.14}

There are eight $h$-fragments of degree~$14$:
\roster
\item\label{h=14,g=1}
  $\GRAPH\AA[3](1;2;3;4) \AA[3](1;2;4;3) \A[2](1,3;2,4) $;
\item\label{h=14,g=2}
  $\GRAPH\AA[3](1;2;3;4) \AA[3](1;2;3;4) \A[2](1,3;2,4) $;
\item\label{h=14,g=3}
  $\GRAPH\AA[4](1;2;3;4;5) \A[3](3,5;2;1,4) \A[1](2,4,5) \(1x3)$;
\item\label{h=14,g=4}
  $\GRAPH\AA[4](1;2;3;4;5) \A[2](1,4;2,5) \A[2](2,4;3,5) \(1x3)$;
\item\label{h=14,g=5}
  $\GRAPH\AA[4](1;2;3;4;5) \A[3](3,4;2;1,5) \A[1](2,4,5) \(1x3)$;
\item\label{h=14,g=6}
  $\GRAPH\AA[4](1;2;3;4;5) \A[2](1,4;2,3) \A[1](1,3,5) \A[1](2,4,5) $;
\item\label{h=14,g=7}
  $\GRAPH\AA[4](1;2;3;4;5) \A[2](1,2;3,4) \A[1](1,3,5) \A[1](2,4,5) $;
\item\label{h=14,g=8}
  $\GRAPH\AA[5](1;2;3;4;5;6) \A[1](1,3,5) \A[1](2,4,6) \(1x4)\(2x5)\(3x6)$,
\endroster
see also \autoref{tab.h=14}.
\table
\caption{$h$-configurations of degree~$14$ \noaux{(see \autoref{conv.tab4})}}\label{tab.h=14}
\vspace{-2\bigskipamount}
\[*
\minitab{\ \ }
\HHH&(r,g,s)&\multispan{21}\ \hss$h$-configurations\hss\\\trule
\iref{h=14,g=1}&(13,4,4)  &2 &1 &\.&\.&\.&\.&\.&\.&\.&\.&\.&\.&\.&\.&\.&\.&\.&\.&\.&\.&\.\\
\iref{h=14,g=2}&(13,4,16) &1 &\.&3 &1 &\.&\.&\.&\.&\.&\.&\.&\.&\.&\.&\.&\.&\.&\.&\.&\.&\.\\
\iref{h=14,g=3}&(14,5,4)  &\.&\.&\.&\.&8 &3 &2 &2 &2 &1 &\.&\.&\.&\.&\.&\.&\.&\.&\.&\.&\.\\
\iref{h=14,g=4}&(14,5,8)  &\.&\.&\.&\.&\.&\.&\.&\.&\.&\.&3 &2 &1 &\.&\.&\.&\.&\.&\.&\.&\.\\
\iref{h=14,g=5}&(14,5,8)  &\.&\.&\.&\.&4 &\.&2 &\.&\.&\.&\.&\.&\.&2 &1 &\.&\.&\.&\.&\.&\.\\
\iref{h=14,g=6}&(14,5,12) &\.&\.&\.&\.&\.&\.&\.&\.&\.&\.&\.&\.&\.&\.&\.&4 &2 &1 &\.&\.&\.\\
\iref{h=14,g=7}&(14,5,14) &\.&\.&\.&\.&\.&\.&\.&\.&\.&\.&\.&\.&\.&\.&\.&\.&\.&\.&1 &\.&\.\\
\iref{h=14,g=8}&(14,6,336)&\.&\.&\.&\.&\.&\.&\.&\.&\.&\.&\.&\.&\.&\.&\.&\.&\.&\.&\.&2 &1 \\
\endminitab
\]
\endtable
The only configuration with twelve $h$-fragments (\viz. eight copies of~$\HH_3$ and
four copies of $\HH_5$) is a maximal $\tD_4$-graph with $26$ vertices.

\subsection{The case $h\sp2=16$}\label{s.16}
There are eight $h$-fragments of degree~$16$:
\roster
\item\label{h=16,g=1}
  $\GRAPH\AA[3](1;2;3;4) \AA[3](1;2;3;4) \AA[3](1;2;4;3) $;
\item\label{h=16,g=2}
  $\GRAPH\AA[3](1;2;3;4) \AA[3](1;2;4;3) \AA[3](1;3;2;4) $;
\item\label{h=16,g=3}
  $\GRAPH\AA[3](1;2;3;4) \AA[3](1;2;3;4) \AA[3](1;2;3;4) $;
\item\label{h=16,g=4}
  $\GRAPH\AA[4](1;2;3;4;5) \AA[4](1;2;3;5;4) \A[1](2,4,5) \(1x3)$;
\item\label{h=16,g=5}
  $\GRAPH\AA[4](1;2;3;4;5) \AA[4](1;2;3;4;5) \A[1](2,4,5) \(1x3)$;
\item\label{h=16,g=6}
  $\GRAPH\AA[4](1;2;3;4;5) \A[3](1,2;4;3,5) \A[3](1,3;4;2,5) $;
\item\label{h=16,g=7}
  $\GRAPH\AA[4](1;2;3;4;5) \A[3](1,4;2;3,5) \A[3](1,3;5;2,4) $;
\item\label{h=16,g=8}
  $\GRAPH\AA[5](1;2;3;4;5;6) \A[2](1,3;2,6) \A[2](3,5;4,6) \(1x4)\(2x5)$,
\endroster
see also \autoref{tab.h=16}.
The configuration with $24$ copies of $\HH_8$ is the only
line-maximizing graph $\Delta'_{32}$, see~\cite{degt:lines}.
Each line is contained in twelve $h$-fragments.
The
$K3$-surface is singular and Kummer: the $32$ lines split
into
two disjoint $16$-tuples of pairwise skew ones, constituting
an abstract configuration $(16_4,16_4)$.

\table
\caption{$h$-configurations of degree~$16$ \noaux{(see \autoref{conv.tab4})}}\label{tab.h=16}
\vspace{-2\bigskipamount}
\[*
\minitab\quad
\HHH&(r,g,s)&\multispan{12}\ \hss$h$-configurations\hss\\\trule
\iref{h=16,g=1}&(14,4,8) &4 &1 &\.&\.&\.&\.&\.&\.&\.&\.&\.&\.\\
\iref{h=16,g=2}&(14,4,24)&\.&\.&1 &\.&\.&\.&\.&\.&\.&\.&\.&\.\\
\iref{h=16,g=3}&(14,4,48)&\.&\.&\.&4 &1 &\.&\.&\.&\.&\.&\.&\.\\
\iref{h=16,g=4}&(14,5,12)&\.&\.&\.&\.&\.&1 &\.&\.&\.&\.&\.&\.\\
\iref{h=16,g=5}&(14,5,12)&\.&\.&\.&\.&\.&\.&1 &\.&\.&\.&\.&\.\\
\iref{h=16,g=6}&(16,5,4) &\.&\.&\.&\.&\.&\.&\.&1 &\.&\.&\.&\.\\
\iref{h=16,g=7}&(16,5,6) &\.&\.&\.&\.&\.&\.&\.&\.&1 &\.&\.&\.\\
\iref{h=16,g=8}&(16,6,96)&\.&\.&\.&\.&\.&\.&\.&\.&\.&24&2 &1 \\
\configurations          &1 &2 &5 &1 &1 &2 &2 &6 &4 &1 &1 &7\\
\endminitab
\]
\endtable

\subsection{The case $h\sp2=18$}\label{s.18}
There are five $h$-fragments:
\roster
\item\label{h=18,g=1}
  $\GRAPH\AA[4](1;2;3;4;5) \AA[4](1;4;2;5;3) \A[3](1,2;4;3,5) $;
\item\label{h=18,g=2}
  $\GRAPH\AA[4](1;2;3;4;5) \AA[4](1;2;5;3;4) \A[3](2,3;5;1,4) $;
\item\label{h=18,g=3}
  $\GRAPH\AA[4](1;2;3;4;5) \AA[4](1;2;5;3;4) \A[3](1,3;5;2,4) $;
\item\label{h=18,g=4}
  $\GRAPH\AA[5](1;2;3;4;5;6) \AA[5](1;2;6;4;5;3) \(1x4)\(2x5)\(3x6)$;
\item\label{h=18,g=5}
  $\GRAPH\AA[5](1;2;3;4;5;6) \AA[5](1;2;3;4;5;6) \(1x4)\(2x5)\(3x6)$,
\endroster
see also \autoref{tab.h=18}, left.
All $h$-fragments are geometric, but not necessarily maximal. The
configuration with three copies of $\HH_5$ is a $\tD_4$-graph with $24$
vertices.

\table
\caption{$h$-configurations of degree~$18$ and $20$ \noaux{(see \autoref{conv.tab4})}}\label{tab.h=18}
\vspace{-2\bigskipamount}
\[*
\minitab\quad
\HHH&(r,g,s)&\multispan6\ \hss$h$-configurations\hss\\\trule
\iref{h=18,g=1}&(16,5,8) &1 &\.&\.&\.&\.&\.\\
\iref{h=18,g=2}&(16,5,8) &\.&1 &\.&\.&\.&\.\\
\iref{h=18,g=3}&(17,5,4) &\.&\.&1 &\.&\.&\.\\
\iref{h=18,g=4}&(17,6,8) &\.&\.&\.&1 &\.&\.\\
\iref{h=18,g=5}&(17,6,24)&\.&\.&\.&\.&3 &1 \\
\configurations          &1 &1 &4 &3 &1 &1 \\
\endminitab
\qquad\quad
\minitab\quad
\HHH&(r,g,s)&\multispan3\ \hss\text{$h$-configs}\hss\\\trule
\iref{h=20,g=1}&(16,6,240)&1&\.&\.\\
\iref{h=20,g=2}&(17,6,48)&\.&1&\.\\
\iref{h=20,g=3}&(18,5,20)&\.&\.&4\\
\configurations          &1 &1 &1 \\
\\\\
\endminitab
\]
\endtable

\subsection{The case $h\sp2=20$}\label{s.20}
There are three $h$-fragments:
\roster
\item\label{h=20,g=1}
  $\GRAPH\AA[5](1;2;3;4;5;6) \AA[5](1;2;3;4;5;6) \A[1](1,3,5) \A[1](2,4,6) $;
\item\label{h=20,g=2}
  $\GRAPH\AA[5](1;2;3;4;5;6) \AA[5](1;3;2;4;6;5) \A[2](2,5;3,6) \(1x4)$;
\item\label{h=20,g=3}
  $\GRAPH\AA[4](1;2;3;4;5) \AA[4](1;2;3;4;5) \AA[4](1;4;2;5;3) $,
\endroster
see also \autoref{tab.h=18}, right.
The graphs $\HH_1$, $\HH_2$ are geometric and maximal, and the graph with
four copies of~$\HH_3$ is $\Phi_{25}$ in~\cite{degt:lines}: it is one of the
three line-maximizing configurations. The other two,
$\Delta'_{25}$ and $\Lambda_{25}$, have no $h$-fragments.

\subsection{The case $h\sp2=22$}\label{s.22}
For each of the remaining values $h^2=22$, $24$, $28$, there is a single
$h$-fragment $\HH_1$, which is contained in a single configuration,
\viz. $\HH_1$ itself,
which is maximal.
For $h^2=22$, the graph $\HH_1$ is encoded as follows:
\roster
\item\label{h=22,g=1}
  $\GRAPH\AA[5](1;2;3;4;5;6) \AA[5](1;4;5;2;6;3) \A[4](3,5;1;2;4,6) $.
\endroster
One has $\rank\Z\HH_1=19$, $\girth\HH_1=6$, and
$\Aut\HH_1\simeq\DG8$.

\subsection{The case $h\sp2=24$}\label{s.24}
The only $h$-fragment~$\HH_1$ in degree~$24$ is
\roster
\item\label{h=24,g=1}
  $\GRAPH\AA[5](1;2;3;4;5;6) \AA[5](1;2;5;6;3;4) \AA[5](1;4;5;2;3;6) $.
\endroster
One has $\rank\Z\HH_1=18$, $\girth\HH_1=6$, and
$\Aut\HH_1\simeq\SG4\times\SG3$.
This graph is $\Lambda_{24}^{\mathrm{A}}$, \ie, the only $3$-regular graph in
\cite[Addendum~3.7]{degt:lines}. It is one of the several line-maximizing
graphs in degree~$24$.

\subsection{The case $h\sp2=28$}\label{s.28}
The only $h$-fragment~$\HH_1$ in degree~$28$ is
\roster
\item\label{h=28,g=1}
  $\GRAPH\DD[5](1,2;3,4;5,6;7,8) \DD[5](1,5;3,7;2,8;4,6) \AA[7](1;7;6;2;3;5;8;4) $.
\endroster
One has $\rank\Z\HH_1=20$, $\girth\HH_1=7$ (the only $h$-fragment~$\HH$ with $\girth\HH>6$),
and
$\Aut\HH_1\simeq\PSL(3,\Bbb{F}_2)\rtimes(\Z/2)$ is of order $336$.
This graph is $\Lambda_{28}$ in~\cite{degt:lines}.
In other words, in degree~$28$, in the only line-maximizing configuration $\Lambda_{28}$
all lines lie in a hyperplane. Simple as it is, this fact was not observed
before.

\section{Hyperelliptic models}\label{S.hyperelliptic}

As a toy example, we discuss briefly smooth hyperelliptic models of
$K3$-surfaces, see \cite[\S\,5.5]{degt:lines}.
Hyperelliptic $h$-fragments are no longer simple graphs, and
they exist only in degrees $h^2\le8$, see \autoref{fig.hyper}.
We summarize the bounds as follows:
\[*
\minitab\quad
h^2=2d\:\ & 2&  4& 6& 8 \\
\max\Hnum\:\ &72&144&36&56
\endminitab
\]

\figure
\def\gap{\ }
\def\Gap{\qquad\quad}
\[*
h^2=2\:\gap
\PIC(0,20)(0,0)
\(0,0){\line(0,1){20}}
\(-1.5,0){\line(0,1){20}}
\(1.5,0){\line(0,1){20}}
\(0,0)\b
\(0,20)\b
\endPIC
\Gap
h^2=4\:\gap
\PIC(15,20)(0,0)
\(-1,0){\line(0,1){20}}
\(1,0){\line(0,1){20}}
\(14,0){\line(0,1){20}}
\(16,0){\line(0,1){20}}
\(0,0){\line(1,0){15}}
\(0,20){\line(1,0){15}}
\(0,0)\b
\(0,20)\b
\(15,0)\b
\(15,20)\b
\endPIC
\Gap
h^2=6\:\gap
\PIC(20,20)(0,0)
\(-1,0){\line(0,1){20}}
\(1,0){\line(0,1){20}}
\(10,0){\line(0,1){20}}
\(19,0){\line(0,1){20}}
\(21,0){\line(0,1){20}}
\(0,0){\line(1,0){20}}
\(0,20){\line(1,0){20}}
\(0,0)\b
\(0,20)\b
\(10,0)\b
\(10,0)\c
\(10,20)\b
\(10,20)\c
\(20,0)\b
\(20,20)\b
\endPIC
\Gap
h^2=8\:\gap
\PIC(40,20)(0,0)
\(-1,0){\line(0,1){20}}
\(1,0){\line(0,1){20}}
\(0,0){\line(1,0){20}}
\(0,20){\line(1,0){20}}
\(10,0){\line(3,2){30}}
\(10,20){\line(3,-2){30}}
\(20,21){\line(1,0){20}}
\(20,19){\line(1,0){20}}
\(20,1){\line(1,0){20}}
\(20,-1){\line(1,0){20}}
\(0,0)\b
\(0,20)\b
\(10,0)\b
\(10,0)\c
\(10,20)\b
\(10,20)\c
\(20,0)\b
\(20,20)\b
\(40,0)\b
\(40,20)\b
\endPIC
\]
\caption{Hyperelliptic $h$-fragments}\label{fig.hyper}
\endfigure

\subsection{Double planes}\label{s.2}
Any degree~$2$ model is hyperelliptic: $X\to\Cp2$ is the double plane
ramified over a smooth sextic curve $C\subset\Cp2$. The projection
establishes a two-to-one correspondence between lines on~$X$ and tritangents
to~$C$ and, essentially by the definition, an $h$-fragment is the pair of
lines over a tritangent. That is, we always have
$\Hnum(X)=\frac12\ls|\Fn X|$. According to~\cite[Addendum~1.2]{degt:sextics},
this number takes values $\{0..66,72\}$ except, possibly, $61$.

\subsection{Double quadrics}\label{s.2x2}
A hyperelliptic quartic $X\to\Cp1\times\Cp1\into\Cp2$ is a double quadric
ramified over a smooth curve $C\subset\Cp1\times\Cp1$ of bi-degree $(4,4)$.
Lines on~$X$ are in a two-to-one correspondence with the generatrices of the
quadric bitangent to~$C$, and an $h$-fragment is the pull-back of
a pair of intersecting bitangents. It follows that $\Hnum=n_1n_2$, where
$n_i$ is the number of bitangents in the $i$-th ruling of the quadric,
$i=1,2$. We have $n_i\le12$, and found in \cite[Corollary~1.4]{degt:singular.K3}
is an example of a quartic with $n_1=n_2=12$. Hence, the sharp upper bound is
$\Hnum(X)\le144$.

\subsection{Sextics and octics}\label{s.hyper}
A hyperelliptic $K3$-surface $X\to Y\into\Cp{d+1}$ of degree $2d=6$ or~$8$ is
a double scroll ramified over a certain smooth curve $C\subset Y$. (We do not
consider the Veronese octic $X\to\Cp2\into\Cp5$ as it has no lines.) The
structure of the graph $\Fn X$ is described in \cite[\S\,5.5]{degt:lines}:
there are
\roster*
\item
a pair of lines $\ell_1,\ell_2$ over the exceptional section $E\subset Y$
(the circled vertices in \autoref{fig.hyper}), so
that $\ell_1\cdot\ell_2=1$ if $2d=6$ and $\ell_1\cdot\ell_2=1$ if $2d=8$, and
\item
a certain number $n$ of pairwise disjoint $\tA_1$-fragments
$\{m_{i1},m_{i2}\}$, one over each bitangent to~$C$, so that
$m_{i1}\cdot m_{i2}=2$ and $m_{ir}\cdot\ell_r=\Gd_{rs}$.
\endroster
We have $n\le9$ if $2d=6$ and $n\le8$ if $2d=8$.

\remark\label{rem.hyper}
In~\cite{degt:lines} it is stated that $n\le10$ if $2d=6$ and $n\le12$ if
$2d=8$. These values can be attained only if the two lines $\ell_1,\ell_2$
over~$E$ are not present. If they are present, the graph is $1$-admissible
but not geometric for $2d=6$ and $10\le n\le12$, and it is not even
hyperbolic for $2d=6$ and $n\ge13$ or $2d=8$ and $n\ge9$.
\endremark

A split hyperplane section on~$X$ is the pull-back of a hyperplane section
of~$Y$ split into~$E$ and $(d-1)$ bitangents.
Since there is a unique $d$-plane in $\Cp{d+1}$ through the line~$E$ and any
$(d-1)$ fibers of the ruling,
the possible counts are the binomial coefficients $C(m,d-1)$, $m\le n$, \ie,
\[*
\Hnum(X)\in
\begin{cases}
  \{ 0, 1, 3, 6, 10, 15, 21, 28, 36\}, & \text{if $2d=6$}, \\
  \{ 0, 1, 4, 10, 20, 35, 56\}, & \text{if $2d=8$}.
\end{cases}
\]

{
\let\.\DOTaccent
\def\cprime{$'$}
\bibliographystyle{amsplain}
\bibliography{degt}
}

\end{document}